\documentclass[12pt]{article}

\usepackage[english]{babel}
\usepackage{amssymb}
\usepackage{amsmath,xypic,epsf,mathrsfs,theorem}

\usepackage[all]{xy}
\newdir{ (}{{}*!/-5pt/\dir^{(}}
\newdir{ )}{{}*!/-5pt/\dir_{(}}

\newtheorem{theorem}{Theorem}[subsection]
\newtheorem{lemma}[theorem]{Lemma}
\newtheorem{corollary}[theorem]{Corollary}

\newtheorem{proposition}[theorem]{Proposition}

{\theorembodyfont{\upshape}
\newtheorem{definition}[theorem]{Definition}

\newtheorem{question}[theorem]{Question}
\newtheorem{problem}[theorem]{Problem}}
\numberwithin{equation}{section}
\numberwithin{theorem}{section}

\newcommand{\Span}{\mathrm{span}}

\newcommand{\cM}{{\cal M}}
\newcommand{\ep}{{\varepsilon}}

\newcommand{\cZ}{{\cal Z}}
\newcommand{\cQ}{{\cal Q}}
\newcommand{\C}{{\mathbb C}}
\newcommand{\Z}{{\mathbb Z}}
\newcommand{\N}{{\mathbb N}}

\newcommand{\Q}{{\mathbb Q}}

\newcommand{\R}{{\mathbb R}}
\newcommand{\Cs}{{C$^*$-al\-ge\-bra}}

\newcommand{\sCs}{{sub-C$^*$-al\-ge\-bra}}

\newcommand{\sh}{{$^*$-ho\-mo\-mor\-phism}}

\newenvironment{proof}[1][Proof.]
{\begin{trivlist}\item[]\textbf{#1} }
{\hbox{}\nobreak\hfill\quad\hbox{$\square$}\end{trivlist}}

\begin{document}

\title{Perturbation of Hausdorff moment sequences, and an
application to the theory of \Cs{}s of real rank zero}

\author{George A.\ Elliott and Mikael R\o rdam}
\date{}
\maketitle

\begin{abstract}
\noindent We investigate the class of unital \Cs s that
admit a unital embedding into \emph{every} unital \Cs{} of real rank
zero, that has no finite-dimensional quotients. We refer to a \Cs{} in
this class as an \emph{initial object}. We show that there are many initial
objects, including for example some unital, simple, infinite-dimensional
AF-algebras, the Jiang-Su algebra, and the GICAR-algebra. 

That the GICAR-algebra is an initial object follows from an analysis
of Hausdorff moment sequences. It is shown that a dense set of 
Hausdorff moment
sequences  belong to a given dense
subgroup of the real numbers, and hence that the Hausdorff
moment problem can be solved (in a non-trivial way) 
when the moments are required to belong to an arbitrary
simple dimension group (i.e.,  unperforated simple ordered 
group with the Riesz decomposition property).
\end{abstract}

\section{Introduction} 
\noindent The following three questions concerning an arbitrary unital
\Cs{} $A$, that is ``large'' in the sense that it has no finite-dimensional
representation, are open.

\begin{question} \label{q1} Does $A$ contain a simple, unital, 
infinite-dimensional sub-\Cs?
\end{question}

\begin{question}[The Global Glimm Halving Problem] \label{q2} 
Does $A$ contain a full\footnote{A subset 
of a \Cs{} is called \emph{full} if it is not contained in any proper 
closed two-sided ideal of the \Cs.} sub-\Cs{} isomorphic to
$C_0((0,1],M_2)$?
\end{question}

\begin{question} \label{q3}
Is there a unital embedding of the Jiang-Su algebra $\cZ$ into $A$?
\end{question}

\noindent The Jiang-Su algebra (see \cite{JiaSu:Z}) is a simple, unital,
infinite-dimensional \Cs, which is $KK$-equivalent to the complex numbers
(and hence at least from a $K$-theoretical point of view could be an initial
object as suggested in Question~\ref{q3}). An affirmative answer to 
Question~\ref{q3} clearly would yield an affirmative answer to 
Question~\ref{q1}. A version of a lemma of
Glimm (see \cite[Proposition~4.10]{KirRor:pi}) confirms Question~\ref{q2}
in the special case that $A$ is simple (and not isomorphic to $\C$);
so Question~\ref{q2} is weaker than Question~\ref{q1}. 

Question~\ref{q2} was raised in \cite[Section~4]{KirRor:pi2} because a 
positive answer will imply that every weakly purely
infinite \Cs{} is automatically purely infinite. 

The Jiang-Su algebra plays a role in the classification program for
amen\-able \Cs s (a role that may well become more important in the future). 
An affirmative answer to Question~\ref{q3} will, besides also answering
the two other questions, shed more light on the Jiang-Su algebra. 
It would for example follow that the Jiang-Su algebra is the (necessarily
unique) unital, simple, separable infinite-dimensional \Cs{} with the
property stipulated in Question~\ref{q3} and with the property (established
in \cite{JiaSu:Z}) that every unital endomorphism can be approximated in
the pointwise-norm topology by inner automorphisms. 

We provide in this paper an affirmative answer to the three questions 
above in the special case in which the target \Cs{} $A$ is required to be
of real rank zero (in addition to being unital and with no finite-dimensional
representations). 

Zhang proved in \cite{Zhang:matricial} that in any unital simple
non-elementary \Cs{} of real rank zero and for any natural number $n$
one can find pairwise orthogonal projections $p_0,p_1, \dots, p_n$ that sum
up to $1$ and satisfy $p_0 \precsim p_1 \sim p_2 \sim \cdots \sim p_n$. In 
other words, one can divide the unit into $n+1$ pieces where $n$ of the
pieces are of the same size, and the last piece is smaller. This result was 
improved in \cite{PerRor:AF} where it was shown that for every natural
number $n$ one can unitally embed $M_n \oplus M_{n+1}$ into any unital \Cs{}
of real rank zero, that has no non-zero representation of dimension $<n$. Thus,
in the terminology of the abstract, $M_n \oplus M_{n+1}$ is an initial
object for every $n$. We shall here extend this result and show that also
the infinite tensor product $P=\bigotimes_{n=1}^\infty M_2 \oplus M_3$ is
an initial object. 

We shall give an algorithm which to an arbitrary unital AF-algebra, that has no 
finite-dimensional representations, assigns a unital \emph{simple} 
infinite-dimensional AF-algebra that embeds unitally into $A$. This leads to the 
existence of a unital infinite-dimensional simple AF-algebra that
unitally embeds into $P$, and hence is an initial object. 
The Jiang-Su algebra was known by Jiang and Su to embed unitally 
into any unital simple non-elementary
AF-algebra, and so is also an initial object.

In Section~\ref{sec:GICAR} we shall show that the Gauge Invariant CAR-algebra
is an initial object. Along the way to this result we shall prove a perturbation
result that may be of independent interest: the set of Hausdorff moment
sequences, with the property that all terms belong to an arbitrary
fixed dense subset of the real numbers, is a dense subset of the
Choquet simplex of all Hausdorff moment sequences. 

We shall show in Section~\ref{sec:4} that a simple, unital,
infinite-dimensional \Cs{} of real rank zero must have 
infinite-dimensional trace simplex if it is an initial object. This leads 
to the open question if one can characterise initial objects among (simple)
unital infinite-dimensional \Cs s of real rank zero (or among simple
AF-algebras). 

We hope that our results will find applications in the future 
study of real rank zero \Cs s; and we hope to have cast some light on
the three fundamental questions raised above.

\section{Initial objects in unital real rank zero \Cs{}s}
\label{sec:initial}

\begin{definition} \label{def:initial}
A unital \Cs{} $A$ will be said in this paper to be an \emph{initial
  object} if it embeds unitally into any unital \Cs{} of real rank
zero which has no non-zero finite-dimensional representations. (Note
that we do not require $A$ to belong to the class of algebras with
these properties.) 
(Also we do not require the embedding to be unique in any way.)
\end{definition}

\noindent
It is clear that the algebra of complex numbers $\C$ is an initial
object, even in the category of all unital \Cs{}s---and that it
is the unique initial object in this larger category. 
It will be shown in Proposition~\ref{prop:initial-i} 
below that the infinite \Cs{} tensor
product $P = \bigotimes_{n=1}^\infty M_2 \oplus M_3$ is also an initial
object in the sense of the present paper. Note that this \Cs{} in fact
belongs to the category we are considering, i.e., unital \Cs{}s
of real rank zero with no non-zero finite-dimensional quotients. 
It follows that a \Cs{} is an
initial object in our sense if and only if it is (isomorphic to) a 
unital sub-\Cs{}
of $P$. We shall use this fact to exhibit a perhaps surprisingly large
number of initial objects, including many simple AF-algebras, the
Jiang-Su algebra, and the GICAR-algebra (the gauge invariant
subalgebra of the CAR-algebra).  

Let us begin by recalling the following standard fact.

\begin{lemma} \label{lm:3}
Let $A$ be a unital \Cs{} and let $F$ be a unital finite-dim\-en\-si\-onal 
sub-\Cs{} of $A$. Let $g_1, \dots, g_n$ denote the minimal (non-zero) 
central projections
in $F$ and let $e_1, \dots, e_n$ be minimal (non-zero) projections in 
$Fg_1, \dots, Fg_n$,
respectively. 

The map consisting of multiplying by the sum $e_1 + \cdots + e_n$ is
an iso\-mor\-phism from the relative commutant $A \cap F'$ of $F$ in $A$ onto the
\sCs{} $e_1Ae_1 \oplus e_2Ae_2 \oplus \cdots \oplus e_nAe_n$ of $A$. Moreover, if $B$
is another unital \Cs{} and $\rho_j \colon B \to e_jAe_j$ are unital
\sh s, then there is a unique  unital \sh{} $\rho \colon B \to A\cap F'$
such that $\rho(b)e_j = e_j\rho(b) = \rho_j(b)$ for all $b \in B$ and
all   $j=1,\dots, n$.
\end{lemma}

\begin{proposition}  \label{prop:initial-i}
The \Cs{} $P=\bigotimes_{n=1}^\infty M_2 \oplus M_3$ is an initial object
(in the sense of Definition~\ref{def:initial}).
\end{proposition}

\begin{proof}
Let $A$ be a unital \Cs{} of real rank zero with no non-zero finite-dimensional
representations. We must find a unital embedding of $P$ into $A$.

Set $\bigotimes_{j=1}^n M_2 \oplus M_3$ = $P_n$, so that $P_{n+1} = P_n \otimes
(M_2 \oplus M_3)$. Let us construct for each $n$ a 
  unital embedding $\varphi_n
\colon P_n \to A$ in such a way that $\varphi_{n+1}(x \otimes 1) =
\varphi_n(x)$ for each $x \in P_n$. This will yield a unital embedding
of $P$ into $A$ as desired.
In order to be able construct these maps inductively we must
require in addition that they be
full.\footnote{By a full
  \sh{} we mean a \sh{} that maps each non-zero element to a full
  element in the codomain algebra. (A full element is one not
  belonging to any proper closed two-sided ideal.)}

For each full projection $e$ in $A$ there is a full unital
embedding\break 
\noindent $\psi \colon M_2 \oplus M_3 \to eAe$. Indeed, $eAe$ cannot
have any non-zero finite-dimensional representation since any such
represen\-tation would extend to a
finite-dim\-en\-si\-onal representation of
$A$ (on a larger Hilbert space).
Hence by \cite[Proposition~5.3]{PerRor:AF}, there
is a unital \sh{} from $M_5 \oplus M_7$ into $eAe$.
Composing this
with a full unital embedding  $M_2 \oplus M_3 \to M_5 \oplus M_7$
yields the desired full embedding $\psi$.

The preceding argument  shows that there is a full unital embedding $\varphi_1
\colon P_1 = M_2 \oplus M_3 \to A$. Suppose that $n \ge 1$ and that
maps $\varphi_1, \varphi_2, \dots, \varphi_n$ have been found 
with the desired properties. 

Choose minimal projections 
$f_1,f_2, \dots, f_{2^n}$
in $P_n$,  one in each minimal non-zero direct  summand, and set 
$\varphi_n(f_j)=e_j $. Each $e_j$ is then a full projection in
$A$. Choose  a full unital embedding $\rho_j \colon M_2
\oplus M_3 \to e_jAe_j$ for each $j$, and note that by Lemma~\ref{lm:3} 
there exists a 
a unital \sh{} $\rho \colon M_2 \oplus M_3 \to A \cap \varphi_n(P_n)'$
such that $\rho(b)e_j = e_j\rho(b) = \rho_j(b)$ for all $b \in
M_2 \oplus M_3$ and all $j$. There is now a unique \sh{} 
$\varphi_{n+1} \colon P_{n+1} = P_n \otimes (M_2
\oplus M_3) \to A$ with the property
that  $\varphi_{n+1}(a
\otimes b) = \varphi_n(a)\rho(b)$ for $a \in P_n$ and $b \in M_2
\oplus M_3$. To show that $\varphi_{n+1}$ is full it suffices to check
that $\varphi_{n+1}(f_j \otimes b)$ is full in $A$ for all $j$ and for all
non-zero $b$ in $M_2 \oplus M_3$;  this follows from the identity
$\varphi_{n+1}(f_n \otimes b) = \varphi_n(f_j)\rho(b) =
e_j\rho(b)=\rho_j(b)$ and the fact that $\rho_j$ is full. 
\end{proof}

\begin{corollary} \label{cor:initial}
Let $A$ be a unital \Cs{} of real rank zero. The following three
conditions are equivalent.
\begin{enumerate}
\item $A$ has no non-zero finite-dimensional representations.
\item There is a unital embedding of $\bigotimes_{n=1}^\infty M_2 \oplus
  M_3$ into $A$.
\item There is a unital embedding of each initial object\footnote{The 
list of initial objects includes some simple unital infinite-dimensional 
AF-algebras and the Jiang-Su algebra $\cZ$ as shown in Section 2.} 
into $A$.
\end{enumerate}
\end{corollary}
 
\begin{proof} (i) $\Rightarrow$ (iii) is true by Definition~\ref{def:initial}. 
(iii) $\Rightarrow$ (ii) follows from Proposition~\ref{prop:initial-i}. (ii)
$\Rightarrow$ (i) holds because any finite-dimensional representation
of $A$ would restrict to a finite-dimensional representation of
$\bigotimes_{n=1}^\infty M_2 \oplus M_3$, and no such exists. 
\end{proof}

\noindent
As remarked above, a \Cs{} is an initial object if and only 
if it embeds unitally into the \Cs{} $P = \bigotimes_{n=1}^\infty M_2
\oplus M_3$. The ordered $K_0$-group of $P$ can be described as follows. Consider
the Cantor set $X = \prod_{n=1}^\infty \{0,1\}$. Consider the maps $\nu_0, \nu_1
\colon X \to \N_0 \cup \{\infty\}$ that for  each $x \in X$
count the number of $0$s and $1$s, respectively, among the
coordinates of $x$, and note that $\nu_0(x)+\nu_1(x)=\infty$ for every $x
\in X$. For each supernatural number $n$ denote by $\cQ(n)$ denote the set
of rational numbers $p/q$ with $q$ dividing $n$, and consider the subgroup  $G \subseteq
C(X,\R)$ consisting of those functions $g$ for 
which $g(x) \in \cQ(2^{\nu_0(x)}3^{\nu_1(x)})$ for every $x \in X$. Equip
$G$ with the pointwise order, i.e., $g \ge 0$ if $g(x) \ge 0$ for all
$x \in X$. Then $(K_0(P),K_0(P)^+,[1])$ is isomorphic to
$(G,G^+,1)$. Note  in particular that $G$ is a dense subgroup of $C(X, \R)$. 

\section{Simple initial objects} \label{sec:2}

\noindent
We shall show in this section that the class of initial objects,
in the sense of 
the previous section, includes several simple
unital (infinite-dimensional) AF-algebras.

\begin{lemma} \label{lm:i}
The following two conditions are equivalent for any dimension group
$G$.
\begin{enumerate}
\item For each order unit $x$ in $G$ there exists an order unit
  $y$ in $G$ such that $2y \le x$.
\item For each finite set of order units $x_1, \dots, x_k$ in $G$
  and for each set of natural numbers $n_1, \dots, n_k$ there is an
  order unit $y$ in $G$ such that $n_jy \le x_j$ for $j=1,2,
  \dots, k$. 
\end{enumerate} 
\end{lemma}

\begin{proof}
The implication
(i) $\Rightarrow$ (ii) follows from the well-known fact
  (which is also  easy to prove---using the
  Effros-Handelman-Shen theorem)  that if $x_1,x_2, \dots, x_k$ are order
  units in a dimension group
  $G$, then there is an order unit $y_0$ in $G$ such that $y_0 \le
  x_j$ for all $j$. The implication  (ii) $\Rightarrow$ (i) is immediate. 
\end{proof}

\noindent 
A dimension group will be said  to have  the property
(D) if it satisfies the two equivalent conditions of
Lemma~\ref{lm:i}. 

\begin{lemma} \label{lm:ii}
Let $A$ be a unital AF-algebra.  The ordered group $K_0(A)$ has the
property (D) if and only if $A$ has no non-zero finite-dimensional representations.
\end{lemma}

\begin{proof}
Suppose that $A$ has no non-zero finite-dimensional representation, and let $x$ be 
an order unit in $K_0(A)$. Then $x = [e]$ for some full projection
$e$ in $M_n (A)$ for some $n$.
Since any finite-dimensional representation of $eM_n(A)e$ would induce
a finite-dimensional representation of $A$ (on a different Hilbert
space), $eM_n(A)e$
has no non-zero finite-dimensional representation. By
\cite[Proposition 5.3]{PerRor:AF} 
there is a unital \sh{} from $M_2
\oplus M_3$ into $eM_n(A)e$. (Cf.~proof of Proposition~\ref{prop:initial-i} above.)
Let $f = (f_1,f_2)$ be a projection in  $M_2 \oplus M_3$, with 
$f_1$ and $f_2$ one-dimensional, and denote by $\tilde{f} \in eM_n(A)e$
the image of $f$ under the unital \sh{} $M_2 \oplus M_3 \to
eM_n(A)e$. Then $\tilde{f}$ is full in $eM_n(A)e$ (because $f$ is full in $M_2
\oplus M_3$), and  $2[\tilde{f}] \le [e]$, as desired. 

Suppose conversely that $K_0 (A)$ has the property (D).
Condition 2.1~(ii) with $k=1$ implies
immediately that every non-zero representation of $A$ is
infinite-dimensional.
\end{proof}

\noindent
We present below a more direct alternative proof (purely in terms of ordered group
theory) of the first implication of the lemma above. 
Consider a decomposition of $K_0 (A)$ as the
ordered group inductive limit of a sequence of
ordered groups $G_1\to G_2 \to \cdots$ with each $G_i$ isomorphic
to a finite ordered group direct sum of copies of $\Z$, 
and let $x$ be an order unit in $K_0(A)$. Modifying
the inductive limit decomposition of $K_0(A)$, we
may suppose that $x$ is the image of an order unit
$x_1$ in $G_1$, and that the image $x_n$ of $x_1$ in $G_n$
is an order unit for $G_n$ for each $n\ge 2$. Let us
show that for some $n$ the
condition~2.1~(i) holds for $x_n$ in $G_n$---or else, if
not, then $G$ has a non-zero quotient ordered group isomorphic to $\Z$.
If not, then for every $n$ there exists
at least one coordinate of $x_n$ in $G_n$ equal to one,
and the inductive
limit of
the sequence consisting, at the $n$th stage, of 
the largest  quotient of  the ordered group $G_n$
in which every coordinate of $x_n$ is equal to one
is a non-zero quotient of $G$ every prime quotient of which is $\Z$.
As soon as Condition~2.1~(i) holds for $x_n$ in $G_n$, then it holds
for $x$ in $G$. In other words, if $G$ has no non-zero quotient isomorphic to $\Z$,
then it has the property (D).

\begin{proposition} \label{prop:ii}
Let $(G,G^+)$ be a dimension group  with  the property
(D). Denote by  $G^{++}$  the set of all order units in $G$, 
and suppose that $G^{++}\not=$ \O. 
Then $(G,G^{++} \cup \{0\})$ is a simple dimension group. 
\end{proposition}

\begin{proof}
Observe first that $G^{++} + G^+ = G^{++}$. With this fact (and with
the assumption that $G^{++}$ is non-empty) it is straightforward to
check that $(G,G^{++} \cup \{0\})$ is an ordered abelian group. We
proceed to show that it is a dimension group. 
This ordered group is unperforated as $(G, G^+)$ is, and so we need only
show that it has the Riesz decomposition property. Equip $G$ with the two
orderings $\le$ and $\precsim$ given by $x \le y$ if $y-x \in G^+$ and
$x \precsim y$ if $y-x \in  G^{++} \cup \{0\}$. Suppose that $x
\precsim y_1 + y_2$ where $x,y_1,y_2 \in G^{++} \cup \{0\}$. We must
find $x_1, x_2 \in  G^{++} \cup \{0\}$ such that $x = x_1+x_2$ and
$x_j \precsim y_j$, $j=1,2$. It is trivial to find $x_1$ and $x_2$ 
in the cases
that one of $x$, $y_1$, $y_2$, and $y_1+y_2-x$ is zero. Suppose 
that the four elements above are non-zero, in which case by hypothesis  they all are
order units. By hypothesis (and by Lemma~\ref{lm:i}) there is $z \in
G^{++}$ such that 
$$2z \le x, \quad z \le y_1, \quad z \le y_2, \quad 2z \le y_1+y_2-x.$$
Then $x-2z \le (y_1-2z) + (y_2-2z)$. Since $(G,G^+)$ has the Riesz
decomposition property there are $v_1,v_2 \in G^+$ such that 
$$x-2z = v_1 + v_2, \qquad v_1 \le y_1-2z, \qquad v_2 \le y_2-2z.$$
Set $ v_1 + z=x_1$ and $ v_2+z=x_2$. Then $x_1, x_2$ belong to
$G^{++}$, $x = x_1+x_2$, $x_1 \precsim y_1$, and $x_2 \precsim y_2$; 
the latter two inequalities hold because 
$$y_j-x_j \; = \; y_j-v_j-z \; = \; (y_j-v_j-2z)+z \; \in \;
G^+ +G^{++} \; = \; G^{++}.$$   
\end{proof}

\begin{proposition} \label{prop:iii}
Let $A$ be a  unital AF-algebra $A$ with no non-zero finite-dimensional 
representation. There exists a unital sub-C$^*$-algebra $B$ of $A$ which 
is a simple, infinite-dimensional AF-algebra, and for which the inclusion 
mapping $B \to A$ gives rise to 
\begin{enumerate}
\item an isomorphism of simplices $T(A) \to T(B)$, and
\item an isomorphism of groups $K_0(B) \to K_0(A)$ which maps $K_0(B)^+$ 
onto $K_0(A)^{++} \cup \{0\}$, and so in particular,
 $$(K_0(B), K_0(B)^+,[1]) \cong (K_0(A), K_0(A)^{++} \cup \{0\}, [1]).$$ 
\end{enumerate}
If $A$ is an initial object, then so also is $B$.
\end{proposition}

\begin{proof} We derive from Lemma~\ref{lm:ii} that $K_0(A)$ has
property (D), and we then conclude from Proposition~\ref{prop:ii} 
that $K_0(A)$ equipped with the 
positive cone $G^+:=K_0(A)^{++}\cup\{0\}$ is a simple 
dimension group. Let $B_1$ be the simple, unital, infinite-dimensional
AF-algebra with dimension group $(K_0(A),G^+,[1_A])$,
and use the homomorphism theorem for AF-algebras 
(\cite[Proposition~1.3.4~(iii)]{Ror:encyc}), to find a
unital (necessarily injective) \sh{} $\varphi \colon B_1 \to A$ which 
induces the (canonical) homomorphism $K_0(B_1) \to K_0(A)$ that maps
$K_0(B_1)^+$ onto $G^+$ and $[1_{B_1}]$ onto $[1_A]$. Set
$\varphi(B_1)=B$. Then $B$ is a unital sub-\Cs{} of $A$, 
$B$ is isomorphic to $B_1$, and (ii) holds.

The property (i) follows from (ii) and the fact, that we shall prove, that the
state spaces of $(K_0(A),K_0(A)^+,[1_A])$ and 
$(K_0(A),G^+,[1_A])$ coincide. The former space is
contained in the latter because $G^+$ is
contained in $K_0(A)^+$. To show the reverse inclusion take a state
$f$ on $(K_0(A),G^+,[1_A])$
and take $g \in K_0(A)^+$. We must show that $f(g) \ge 0$. 
Use Lemmas~\ref{lm:i} and
\ref{lm:ii} to find for each natural number $n$ an element $v_n$ in
$K_0(A)^{++}$ such that $nv_n \le [1_A]$. Then $nf(v_n) \le 1$, so
$f(v_n) \le 1/n$; and $g+v_n$ belongs to
$K_0(A)^{++}$, so $f(g+v_n) \ge 0$. These two inequalities, that hold for
all $n$, imply that $f(g) \ge 0$. 
\end{proof}

\begin{corollary} \label{cor:initial_ii}
\mbox{}
\begin{enumerate}
\item There is a simple unital infinite-dimensional AF-algebra which
  is an initial object.
\item The Jiang-Su algebra $\cZ$ is an initial object.
\end{enumerate}
\end{corollary}

\begin{proof}
The assertion (i) follows immediately  from
Propositions~\ref{prop:initial-i} and \ref{prop:iii}. 

The assertion (ii) follows from (i) and the fact, proved in \cite{JiaSu:Z}, that the
Jiang-Su algebra $\cZ$ embeds in (actually is tensorially absorbed by) any
unital simple infinite-dimensional AF-algebra.
\end{proof}

\noindent
The corollary above provides an affirmative answer to
Question~\ref{q3} (and hence also to Questions~\ref{q1} and \ref{q2})
of the introduction in the case that the target \Cs{} $A$ is
assumed to be of real rank zero. 

The question of initial objects may perhaps be pertinent in the classification
program, where properties such as approximate divisibility and being able to
absorb the Jiang-Su algebra $\cZ$ are of interest. We remind the reader that a
\Cs{} $A$ is approximately divisible if for each natural number $n$ there is a 
sequence $\varphi_k \colon M_n \oplus M_{n+1} \to \cM(A)$ of unital \sh s
(where $\cM(A)$ denotes the multiplier algebra of $A$) such that $[\varphi_k(x),a] 
\to 0$ for all $a \in A$ and all $x \in M_n \oplus M_{n+1}$. (It turns out
that if $A$ is unital, then we need only find such a sequence of \sh s for
$n=2$.) It is easily seen that a separable \Cs{} $A$ is approximately 
divisible if, and only if, there is a unital \sh{}
\begin{equation} \label{eq:prod}
\prod_{n\in \N} (M_n \oplus M_{n+1}) / \sum_{n \in \N} (M_n \oplus M_{n+1})
\to \cM(A)_\omega \cap A',
\end{equation}
and it follows from \cite[Theorem~7.2.2]{Ror:encyc} and \cite{JiaSu:Z}
that $A$ is $\cZ$-absorbing 
if and only if there is  a unital embedding of $\cZ$ into $\cM(A)_\omega \cap
A'$; here, $\omega$ is any free ultrafilter on $\N$, 
and $\cM(A)$ is identified with a sub-\Cs{} of the ultrapower $\cM(A)_\omega$
(the \Cs{} of bounded sequences in $\cM(A)$, modulo the
ideal of bounded sequences convergent to $0$ along $\omega$).

Toms and Winter recently observed (\cite{TomsWin:Z}) that
any separable approximately divisible \Cs{} is $\cZ$-absorbing, because one
can embed $\cZ$ unitally into the \Cs{} on the left-hand side of
\eqref{eq:prod}. (The latter fact follows from our 
Corollary~\ref{cor:initial_ii}, but it can also be proved directly, as was done
in \cite{TomsWin:Z}.)
In the general case, when $A$ need not be
approximately divisible, it is of interest to decide when $A$ is 
$\cZ$-absorbing, or, equivalently, when one can find a unital embedding of
$\cZ$ into $\cM(A)_\omega \cap A'$.  Here it would  be
extremely useful if one knew that $\cZ$ was an initial object in the
category of all unital \Cs s with no non-zero finite-dimensional representations. 

The proof of Corollary~\ref{cor:initial_ii} yields an explicit---at the level of
the invariant---simple
unital AF-algebra which is an initial object. Indeed, consider the
initial object $P = \bigotimes_{n=1}^\infty M_2 \oplus M_3$, the 
$K_0$-group of which is the dense subset $G$ of $C(X,\R)$ described 
above (after Corollary~\ref{cor:initial}), with the relative order,
where $C(X,\R)$ is equipped with the
standard pointwise ordering. The simple dimension group 
$(G,G^{++} \cup \{0\})$ of Proposition~\ref{prop:ii} is obtained
by again viewing  $G$ as a subgroup of $C(X,\R)$ but this time endowing 
$C(X,\R)$ with the strict  pointwise ordering (in which an element $f
\in C(X,\R)$ is  
positive if $f =0$ or if $f(x) >0$ for all $x \in X$). Any other
simple dimension group which maps onto this may also be used. 

It would of course be nice to have an even more explicit (or natural) 
example of a simple unital infinite-dimensional AF-algebra which is 
an initial object in the sense of this paper.

The trace simplex of the simple unital AF-algebra referred to  above
is the simplex of 
probability measures on the Cantor set. We shall show in Section~\ref{sec:4} that
the trace simplex of an initial object, that has sufficiently many
projections,  must be infinite-dimensional. 
Let us now note that a large class of infinite-dimensional 
Choquet simplices  arise as the trace simplex of an
initial object.

\begin{proposition} \label{prop:initial}
Let $X$ be a metrizable compact Hausdorff space which admits an embedding of the
Cantor set{}.\footnote{An equivalent formulation of this
(rather weak) property is that $X$ has a non-empty closed
subset with no isolated points.}
There exists a simple unital AF-algebra $A$
which is an initial object, such that $T(A)$ is affinely homeomorphic
to the simplex $\cM_1(X)$  of (Borel) probability measures on $X$.  
\end{proposition}

\begin{proof}
By hypothesis $X$ has a closed subset $X_0$ which is (homeomorphic to)
the Cantor set. The dimension group of the known initial object 
$\bigotimes_{n=1}^\infty M_2 \oplus M_3$ is isomorphic in a natural
way to a dense subgroup $G$ of
$C(X_0,\R)$ (equipped with the standard pointwise ordering), with
canonical order unit corresponding to  
the constant function $1_{X_0}$, cf.\ the remark after
Corollary~\ref{cor:initial}. We shall construct below
a countable dense subgroup $H$ of $C(X,\R)$ such that the constant
function $1_X$ belongs to $H$, and such that the restriction
$f|_{X_0}$ belongs to $G$
for every $f \in H$. Equip $H$ with the strict pointwise ordering on $C(X,\R)$
and with the order unit $1_X$.  
Then we have an ordered  group homomorphism
$\varphi \colon H \to G$ given by $\varphi(f) = f|_{X_0}$, which maps
$1_X$ into $1_{X_0}$. It follows that we
may take $A$ to be the unital, simple AF-algebra with invariant
$(H,H^+,1_X)$, as 
by the homomorphism theorem for AF-algebras (cf.~above)
$\varphi$ induces a unital embedding of $A$ into
$\bigotimes_{n=1}^\infty M_2 \oplus M_3$, whence $A$ is an initial object, and the trace
simplex of $A$ is homeomorphic to the state space of $(H,H^+,1_X)$,
which is $\cM_1(X)$.

Let us now pass to the construction of $H$.  Each $g \in G$ extends to $\tilde{g} \in
C(X,\R)$ (we do not make any assumption concerning the mapping $g \mapsto \tilde{g}$). 
Choose a countable dense subgroup $H_0$ of 
$C_0(X\setminus X_0, \R) \subseteq C(X,\R)$, and consider the countable subgroup of
$C(X,\R)$ generated by $H_0$
and the countable set $\{\tilde{g} : g \in G\}$. 
Denote this group, with the relative (strict pointwise) order, by $H$; let us
check that this choice of $H$ fulfils the requirements.
First, $f|_{X_0} \in G$
for every $f \in H$. To 
see that $H$ is dense in $C(X,\R)$, 
let there be given
$f \in C(X,\R)$ and  $\ep
>0$. Choose  $g \in G$ such that $\|f|_{X_0}-g\|_\infty < \ep/2$. Extend
$f|_{X_0}-g$ to a function $f_0 \in C(X,\R)$ with $\|f_0\|_\infty =
\|f|_{X_0}-g\|_\infty < \ep/2$. Note that $f-\tilde{g}-f_0$ belongs to
$C_0(X\! \setminus\! X_0,\R)$. Choose $h_0 \in H_0$ such that
$\|f-\tilde{g}-f_0-h_0\|_\infty < \ep/2$, and consider the function $h=\tilde{g}+h_0 \in
H$. We have  $\|f-h\|_\infty \le \|f-\tilde{g}-f_0-h_0\|_\infty +
\|f_0\|_\infty < \ep$, as desired.
\end{proof}

\section{Hausdorff moments, the GICAR-algebra, and Pascal's triangle} \label{sec:GICAR}

\noindent In this section we shall establish  the following result.

\begin{theorem} \label{thm:GICAR} The GICAR-algebra is an initial
  object (in the sense of Definition~\ref{def:initial}).
\end{theorem}

\noindent We review some of the background material.
Consider the Bratteli diagram given by Pascal's triangle, 
 \setlength{\unitlength}{1mm}
  \newcommand{\linie}{\ar@{-}}     
 \newcommand{\dlinie}{\ar@{=}}
  \newcommand{\entry}[1]{\phantom{^{\textstyle#1}}
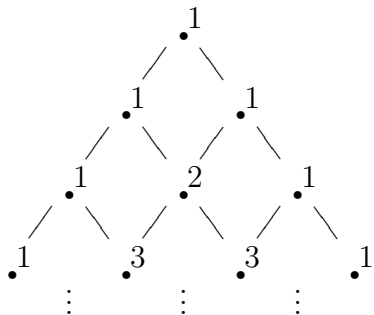
    \begin{picture}(1,0)\put(.5,.5){\circle*{1}}%
    \end{picture}^{\textstyle#1}}
\begin{displaymath}
    \xymatrix@R-1pc@C-2pc{ &&&&{\entry1}\linie[dl]\linie[dr]\\ 
      &&&{\entry1}\linie[dl]\linie[dr]&&{\entry1}\linie[dl]\linie[dr]\\ 
      &&{\entry1}\linie[dl]\linie[dr]&&{\entry2}\linie[dl]\linie[dr]&&
      {\entry1}\linie[dl]\linie[dr]\\ 
      &{\entry1}&&{\entry3}\ar@{}[dll]|{\objectstyle\vdots}
      &&{\entry3}\ar@{}[dll]|{\objectstyle\vdots}&&
      {\entry1}\ar@{}[dll]|{\objectstyle\vdots}\\&&&&&}
  \end{displaymath}
and denote by 
$$\C = B_0 \to B_1 \to  B_2 \to \cdots \to \varinjlim B_n \; (\,=B)$$
the inductive system of finite-dimensional \Cs s associated with
that Bratteli diagram. The \Cs{} $B$ is the GICAR-algebra. (It can
also, more naturally, be realized as the fixed point algebra of the 
CAR-algebra under a certain action of the circle referred to as the
gauge invariant action, cf.\ \cite{Bra:UHF}.)

For each $n \ge 0$ and $0 \le k \le n$, choose a minimal projection 
$e(n, k)$ in the $k$th minimal direct summand of $B_n$. Note that 
$e(0,0)=1_B$ and that $e(n,k)$ is Murray-von Neumann equivalent to 
$e(n+1,k) + e(n+1,k+1)$ in $B_{n+1}$. A trace $\tau$ on $B_n$ is determined by 
its values on  the projections $e(n,k)$, $0 \le k \le n$. 

The group $K_0(B)$ is generated, as an ordered abelian group, by the elements
$[e(n,k)]$, with $n \ge 0$ and $0 \le k \le n$; that is, these elements span
$K_0(B)$ as an abelian group, and the semigroup spanned by the
elements $[e(n,k)]$ is equal to $K_0(B)^+$. Our
generators satisfy the relations
\begin{equation} \label{eq:relation}
[e(n,k)] = [e(n+1,k)] + [e(n+1,k+1)], \qquad n \ge 0, \; \, 0 \le k \le n.
\end{equation}
Moreover, $(K_0(B),K_0(B)^+)$ is the universal ordered abelian group
generated, as an ordered abelian group, by elements $g(n,k)$, $n \ge 0$ and $0
\le k \le n$, with the relations $g(n,k) = g(n+1,k)+g(n+1,k+1)$. 

For brevity we shall set $(K_0(B),K_0(B)^+,[1_B]_0) = (H,H^+,v)$. 

For each abelian (additively written) 
group $G$ and for each sequence $t \colon \N_0 \to G$ 
associate the discrete derivative $t' \colon \N_0 \to G$ given by $t'(k)
= t(k)-t(k+1)$. Denote the $n$th derivative of $t$ by  $t^{(n)}$, and apply
the convention $t^{(0)} = t$. 

We remind the reader of the following classical result. The equivalence of
(i) and (iv) is the solution to the Hausdorff Moment problem (see e.g.\
\cite[Proposition~6.11]{BergChrRes:Harmonic}). The equivalence of (i), (ii),
and (iii) follows from Proposition~\ref{prop:b1} below (with $(G,G^+,u) =
(\R,\R^+,1)$). 

\begin{proposition}[Hausdorff Moments]
\label{prop:b}
The following four conditions are equivalent for any sequence
$t \colon \N_0 \to \R$.
\begin{itemize}
\item[{\rm{(i)}}] $t^{(k)}(n) \ge 0$ for all $n,k \ge 0$.
\item[{\rm{(ii)}}] There is a system, $\{t(n,k)\}_{0 \le k \le n}$,
of positive real numbers (necessarily unique) such that
$$t(n+1,k) + t(n+1,k+1) = t(n,k), \qquad t(n,n) = t(n),$$
for $n \ge 0$ and $0 \le k \le n$. 
\item[{\rm{(iii)}}] There is a (unique) tracial state $\tau$ on the GICAR-algebra 
such that $t(n) = \tau(e(n,n))$ for all $n \ge 0$.   
\item[{\rm{(iv)}}] There is a Borel probability measure $\mu$ on the interval
  $[0,1]$ such that 
$$t(n) = \int_0^1 \lambda^n \, d\mu(\lambda),$$
for all $n \ge 0$. 
\end{itemize} 
\end{proposition}

\vspace{.3cm} \noindent It follows from Proposition~\ref{prop:b1}
below and from (iv) that the coefficients
$t(n,k)$ from (ii) are given by
\begin{equation} \label{eq:t(n,k)}
t(n,k) = t^{(n-k)}(k) = \int_0^1 \lambda^k (1-\lambda)^{n-k} \,  d\mu(\lambda).
\end{equation}

\noindent
A sequence $t = (t(0),t(1), \dots)$ satisfying the condition in
Proposition~\ref{prop:b}~(iv) (or, equivalently, the three other conditions
of Proposition~\ref{prop:b}) is called
a \emph{Hausdorff moment sequence}. Note that $t(0)=1$ in every Hausdorff
moment sequence. Let us denote the set of all moment sequences 
by $\cM$. Note that $\cM$ is a compact convex set and in fact a Choquet 
simplex. For each $n \in \N_0$ let us set
$$\cM_n = \big\{\big(t(0),t(1),t(2),  \dots, t(n)\big) :
\big(t(0),t(1),t(2), \dots \big) \in \cM\big\} \subseteq \R^{n+1},$$
and denote by $\pi_n$ the canonical surjective affine mapping
$\cM_{n+1} \to \cM_n$. 

Let us say that a moment sequence $t=(t(0),t(1),t(2), \dots)$ is 
 \emph{trivial} if the corresponding measure in
 Proposition~\ref{prop:b}~(iv) is supported in $\{0,1\}$, and say that $t$
 is \emph{non-trivial} otherwise. A sequence $t$ is trivial if and
 only if it is a convex combination of the two trivial sequences
 $(1,1,1, \dots)$ and $(1,0,0,\dots)$. It follows from this and (iv)
 above that $t$ is non-trivial if and only if $t(2) < t(1)$. One
 can use Equation~\eqref{eq:t(n,k)} to see that $t$ is non-trivial if
 and only if $t(n,k) \ne 0$ for all $n$ and $k$. 

We seek unital embeddings from the GICAR algebra $B$ into unital AF-algebras
(and into unital \Cs s of real rank zero). At the level of the invariant
we are thus seeking positive unit preserving group homomorphisms from the 
dimension group with distinguished unit $(H,H^+,v)$ associated to the
GICAR algebra into the 
ordered $K_0$-group with distinguished unit of the target algebra; call
this invariant $(G,G^+,u)$. The proposition below rephrases this problem
as that of the existence of a function $g \colon \N_0 \to G$ with certain 
properties. 

\begin{proposition} \label{prop:b1}
Let $(H,H^+,v)$ be as above, and let $(G,G^+,u)$ be an ordered abelian 
group with a distinguished order unit $u$. Let $g \colon \N_0 \to G$ be 
given, and assume that $g(0)=v$.  
The following  conditions are equivalent.
\begin{enumerate}
\item $g^{(k)}(n) \in G^+$ for all $n,k \ge 0$.
\item There is a system, $\{g(n,k)\}_{0 \le k \le n}$,
of elements in $G^+$ (necessarily unique) such that
$$g(n+1,k) + g(n+1,k+1) = g(n,k), \qquad g(n,n) = g(n),$$
for all $n \ge 0$ and $0 \le k \le n$.
\item There is a (unique) homomorphism of ordered groups 
$\varphi \colon H \to G$ with $\varphi(v)=u$ 
such that $\varphi([e(n,n)]) = g(n)$ for all $n \ge 0$.
\end{enumerate} 
If the three conditions above are satisfied, then 
$$\varphi([e(n,k)]) = g(n,k) = g^{(n-k)}(k)$$
for all $n \ge 0$ and $0 \le k \le n$; and the homomorphism $\varphi$
is faithful if and only if $g(n,k)$ is non-zero for all $n \ge 0$ and $0
\le k \le n$.
\end{proposition}

\begin{proof} (i) $\Rightarrow$ (ii). Set $g(n, k) = g^{(n-k)}(k) \in G^+$. 
Then $g(n,n) = g^{(0)}(n) = g(n)$, and 
\begin{eqnarray*}
g(n,k)-g(n+1,k+1) & = & g^{(n-k)}(k) - g^{(n-k)}(k+1) \; = \;
g^{(n-k+1)}(k) \\ & = & g(n+1,k).
\end{eqnarray*}

(ii) $\Rightarrow$ (iii). We noted after Theorem~\ref{thm:GICAR} that 
$H=K_0(B)$ is generated, as an ordered abelian group, by the elements 
$[e(n,k)]$, $n \ge 0$ and $0 \le k \le n$, and that $H$ is the
universal ordered abelian group generated by these elements subject to
the relations \eqref{eq:relation}. Accordingly, by (ii), there
exists a (unique) positive group homomorphism $\varphi \colon H \to G$
with $\varphi([e(n,k)]) = g(n,k)$. Also,
$\varphi(v) = \varphi([e(0,0)]) = g(0,0) = g(0) = u$.

To complete the proof we must show that $\varphi$ is uniquely determined by
its value on the elements $[e(n,n)]$, $n \ge 0$. But this follows from the fact
that the elements $[e(n,k)]$, with $n \ge 0$ and $0 \le k \le n$, 
belong to the subgroup generated by the elements $[e(n,n)]$, for $n \ge 0$, by the
relations \eqref{eq:relation}.

(iii) $\Rightarrow$ (i). This implication follows from the identity 
$\varphi([e(n+k,n)])=g^{(k)}(n)$, that we shall proceed to prove by
induction on $k$. The case $k=0$ is explicitly contained in (iii). 
Assume that the identity has been shown to hold for some $k \ge 0$. Then,
by \eqref{eq:relation},
\begin{eqnarray*}
g^{(k+1)}(n) & = & g^{(k)}(n) - g^{(k)}(n+1) = \varphi\big([e(n+k,n)] -
[e(n+k+1,n+1)]\big) \\
&=& \varphi([e(n+k+1,n)]).
\end{eqnarray*}

To prove the two last claims of the proposition, assume that $g$ satisfies
the three equivalent conditions, and consider the homomorphism of
ordered groups $\varphi
\colon H \to G$ asserted to exist in (iii). It follows from the proofs
of (i) $\Rightarrow$ (ii) and (ii) $\Rightarrow$ (iii) that 
$\varphi([e(n,k)]) = g(n,k) = g^{(n-k)}(k)$. Any non-zero positive
element $h$ of $H$ is a finite (non-empty) sum of elements of the form
$[e(n,k)]$. Thus $\varphi(h)$ is a finite (non-empty) sum of elements
of the form $g(n,k)$. This shows that $\varphi(h)$ is non-zero for all non-zero
positive elements $h$ in $H$ if and only if $g(n,k)$ is non-zero for
all $n$ and $k$.
\end{proof}

\noindent Let us now return to the convex set $\cM$ of Hausdorff 
moment sequences in $\R^+$ and to the truncated finite-dimensional
convex sets $\cM_n$. 

\begin{lemma} \label{lm:c}
$\dim(\cM_n) = n$.
\end{lemma} 

\begin{proof}
The convex set $\cM_n$ is a subset of $\{1\} \times \R^n$, and has
therefore dimension at most $n$. On the other hand, the
points $(1,\lambda,\lambda^2, \dots, \lambda^n)$ belong to $\cM_n$
for each $\lambda \in (0,1)$, and these points span an $n$-dimensional
convex set. 
\end{proof}

\noindent Let $\cM_n^\circ$ denote the relative interior\footnote{The
  relative interior of a finite-dimensional convex set is its interior
  relatively to the affine set it generates.} of
$\cM_n$. By standard theory for finite-dimensional convex sets, see
e.g.\ \cite[Theorem~3.4]{Bro:convex}, $\dim(\cM_n^\circ) = \dim(\cM_n) = n$. 
Note that
$$\cM_1 = \{(1,\lambda) : \lambda \in [0,1] \}, \qquad \cM_1^\circ =
\{(1,\lambda) : \lambda \in (0,1)\}.$$
For $n \ge 2$ we can use Lemma~\ref{lm:c} to conclude that $\cM_n^\circ =
\{1\} \times U_n$ for some open convex subset $U_n$ of $\R^n$.  

\begin{lemma} \label{lm:b}
$\pi_n(\cM_{n+1}^\circ)=\cM_n^\circ$.
\end{lemma}

\begin{proof} This follows from the standard fact from the theory for 
finite-dimensional convex sets (see e.g.\ 
\cite[\S 3 and Exercise 3.3]{Bro:convex}) that the relative interior
of the image of $\pi_n$ is the image under $\pi_n$ of the relative interior
of $\cM_{n+1}$ (combined with the fact that $\pi_n$ is surjective).
\end{proof}

\begin{theorem} \label{thm:a}
Let $G$ be a dense subset of the reals that contains $1$. Then there
is a non-trivial moment sequence $(t_0,t_1,t_2, \dots)$ such that $t_n$
belongs to $G$ for every $n \in \N_0$. 
Furthermore, the moment sequences with all terms belonging
to $G$ constitute a dense\footnote{In the standard pointwise (or
  product) topology.} subset of $\cM$. If $G$ also is a group, and has
infinite rank over $\Q$, then there exists a  moment sequence in $G$
the terms of which are independent over $\Q$.
\end{theorem}

\begin{proof} Let $(s_0, s_1, s_2, \dots)$ be a moment sequence, let
$m$ be a natural number, and let $\ep_1, \ep_2, \dots, \ep_m$ be 
strictly positive real numbers. Since $(s_0,s_1, \dots, s_m)$ belongs
to $\cM_m$, since $\cM_m^\circ$ is dense in $\cM_m$ (cf.\
\cite[Theorem~3.4]{Bro:convex}) and is equal to $\{1\} \times U_m$
for some open subset $U_m$ of $\R^m$, since $1 \in G$, and since $G$ is dense in $\R$, 
we can find $(t_0,t_1, \dots, t_m)$ in $\cM_m^\circ$ such that 
$t_j$ belongs to $G$ for $j=0,1,\dots,m$ and $|t_j-s_j| < \ep_j$ for $j=1, \dots,
m$. 

Let us choose inductively $t_n$, $n > m$, such that $t_n \in G$ and
$(t_0,t_1, \dots, t_n) \in \cM_n^\circ$. Supppose that $n \ge m$ and
that $t_0, t_1, \dots, t_n$ have been found. The set
$$\{s \in \R : (t_0,t_1, \dots, t_n,s) \in \cM_{n+1}^\circ\}$$
is non-empty (by Lemma~\ref{lm:b}) and open 
(because $\cM_{n+1}^\circ = \{1\} \times U_{n+1}$
for some open subset $U_{n+1}$ of $\R^{n+1}$). 
Hence there exists  $t_{n+1} \in G$ such that 
 $(t_0,t_1, \dots, t_{n+1}) \in \cM_{n+1}^\circ$.

The resulting sequence $(t_0, t_1, t_2,\dots)$ in $G$ is
a moment sequence by construction and is close to 
the given moment sequence $(s_0, s_1, s_2, \dots)$.

The inequality $t_2<t_1$ holds because $(t_0,t_1,t_2)$
belongs to the open set $\cM_2^\circ = \{1\} \times U_2$. (Indeed, note that
$t_1 \le t_2$ whenever $(t_0,t_1,t_2)$ belongs to $\cM_2$ and, hence, that
the element $(1,t_1,t_1)$ of $\cM_2$ belongs to the boundary.)

Concerning the desired independence of the terms of the moment
sequence when $G$ is a group, of infinite rank, it will suffice to
choose each $t_n$ in the set
$$G \setminus \Span_\Q\{t_0,t_1, \dots, t_{n-1}\}.$$
This is possible because this set is dense in $\R$ by
the assumption on $G$. 
\end{proof}

\begin{corollary} \label{cor:b}
Let $G$ be a dense subgroup of $\R$ with $1 \in G$. There is a faithful 
homomorphism of ordered groups from the dimension group
$H$ associated with the Pascal triangle to $G$  
(with the order inherited from $\R$)
that maps the canonical  order unit of $H$  to $1$. 
Furthermore, the set of such maps into $G$ is dense in
the set of such maps just into $\R$, in 
the topology of pointwise convergence on $H$.
If $G$ is of infinite rank there 
exists such a map which is injective.
\end{corollary}

\begin{proof}
Propositions~\ref{prop:b} and \ref{prop:b1} give a one-to-one
correspondence between moment sequences $t \colon \N_0 \to G$ and
homomorphisms $\varphi \colon H \to G$ of ordered abelian groups that
map the canonical order unit $v \in H$ into $1 \in G$, such that
$\varphi([e(n,k)]) = t(n,k)$ for all $n \ge 0$ and $0 \le k \le n$. 
If $t$ is non-trivial, then $t(n,k)$ is
non-zero for all $n,k$, whence $\varphi(g) >0$ for every non-zero
positive element $g$ in $H$ (because each such element $g$ is a sum
of elements of the form $[e(n,k)]$). 

A pointwise converging net of moment sequences
corresponds to a pointwise converging net of homomorphisms $H \to
G$. 

The first two claims now follow from Theorem~\ref{thm:a}.

A homomorphism $\varphi \colon H \to G$ is injective if the
restriction of $\varphi$ to the sub-group spanned by $\{[e(n,k)] : k =
0,1, \dots,n\}$ is injective for every $n$. The latter holds, for
a specific $n$, if and only if 
$t(n,0), t(n,1), \dots, t(n,n)$ are independent over $\Q$, or,
equivalently, if and only if $t(0), t(1), \dots, t(n)$ are independent over
$\Q$. (Use the relation in Proposition~\ref{prop:b}~(ii) to see the
second equivalence.) This shows that a moment sequence $t \colon \N_0 \to G$, where
$t(0),t(1),\dots$ are independent over $\Q$, gives rise to an injective
homomorphism $\varphi \colon H \to G$. The existence of such a moment
sequence $t$, under the assumption that $G$ has infinite rank, follows
from Theorem~\ref{thm:a}.
\end{proof}

\begin{lemma}  \label{lm:d}
With  $X$ the Cantor set, let $f_1, \dots, f_n \colon X \to \R$ be
continuous functions, and let $U \subseteq \R^{n+1}$ be an open subset
such that
$$\{s \in \R : (f_1(x),f_2(x), \dots, f_n(x),s) \in U\}$$
is non-empty for every $x \in X$. It follows that there exists a
continuous function $f_{n+1} \colon X \to \R$ such that
$$(f_1(x),f_2(x), \dots, f_n(x), f_{n+1}(x)) \in U$$
for all $x \in X$.
\end{lemma}

\begin{proof}
For each $s \in \R$ consider the set $V_s$  of those $x \in X$ for which 
$(f_1(x),f_2(x), \dots, f_n(x),s)$ belongs to $U$. Then $(V_s)_{s \in \R}$ is
an open cover of $X$, and so by compactness, $X$ 
has a finite subcover $V_{s_1}, V_{s_2}, \dots, V_{s_k}$. Because $X$
is totally disconnected there are clopen subsets $W_j \subseteq V_{s_j}$
such that $W_1, W_2, \dots, W_k$ partition $X$. The function $f_{n+1} =
\sum_{j=1}^k s_j 1_{W_j}$ is as desired.
\end{proof}

\begin{proposition} \label{prop:c}
With $X$  the Cantor set,  let $G$ be a norm-dense subset of
$C(X,[0,1])$ that contains the constant function $1$. 
There exists a sequence $(g_0,g_1,g_2, \dots)$ in $G$ 
such that $(g_0(x),g_1(x),g_2(x), \dots)$ is a
non-trivial moment sequence for every $x \in X$. 
\end{proposition}

\begin{proof}
Choose $g_0,g_1,\dots$ in $G$ inductively such that
$(g_0(x),g_1(x), \dots, g_n(x))$ belongs to $\cM_n^\circ$ for every $x
\in X$. Begin by choosing $g_0$ to be the constant function 1 (as it
must be). 
Suppose that $n \ge 0$ and that $g_0,g_1, \dots, g_n$ as above have 
been found. As observed earlier, $\cM_{n+1}^\circ = \{1\} \times
U_{n+1}$ for some open subset $U_{n+1}$ of $\R^{n+1}$. The set 
\begin{eqnarray*} && \{s \in \R : (g_1(x), \dots, g_n(x),s) \in U_{n+1} \} \\
&=& \{s \in \R : (g_0(x),g_1(x), \dots, g_n(x),s) \in \cM_{n+1}^\circ \}
\end{eqnarray*}
is non-empty for each $x \in X$ (by Lemma~\ref{lm:b}), 
and so we can use Lemma~\ref{lm:d} to find a continuous
function $f \colon X \to \R$ such that $(g_1(x),\dots, g_n(x),f(x))$
belongs to $U_{n+1}$ for all $x \in X$. By compactness of $X$, continuity of
the functions $g_1, \dots, g_n,f$, and because $U_{n+1}$ is open, there
exists $\delta >0$ such that $(g_1(x),\dots, g_n(x),h(x))$
belongs to $U_{n+1}$ for
all $x \in X$ whenever $\|f-h\|_\infty < \delta$. As $G$ is dense in
$C(X,\R)$ we can find $g_{n+1} \in G$ with $\|f-g_{n+1}\|_\infty <
\delta$, and this function has the desired properties.

As in the proof of Proposition~\ref{prop:b}, since 
$(g_0(x),g_1(x),g_2(x))$ belongs to $\cM_2^\circ$, we get $g_2(x) < g_1(x)$, 
which in turns implies that the moment sequence $(g_0(x),g_1(x),g_2(x), \dots)$ is 
non-trivial for every $x \in X$. 
\end{proof}

\begin{proposition} \label{prop:d}
With $X$  the Cantor set, let $G$ be a norm-dense subgroup of
$C(X,\R)$ that contains the constant function $1$. 
There exists a faithful homomorphism of ordered groups from the 
dimension group $H$ associated with the Pascal triangle to $G$ 
(with the strict pointwise order)
that  takes the distinguished order unit $v$ of $H$ into the 
constant function $1$. 
\end{proposition}

\begin{proof}
Choose a sequence $g_0,g_1,g_2,\dots$ in $G$ as specified in
Proposition~\ref{prop:c}, and consider the (unique) system
$\{g(n,k)\}_{0 \le k \le n}$ in $G$ such that 
$$g(n+1,k) + g(n+1,k+1) = g(n,k), \qquad g(n,n) = g_n$$ 
for $n \ge 0$ and $0 \le k \le n$. Use Proposition~\ref{prop:b} and
the non-triviality of the moment sequence $(g_0(x),g_1(x),g_2(x), \dots)$ 
to conclude that $g(n,k)(x) > 0$ for all $x \in X$. Hence, by 
Proposition~\ref{prop:b1}, there exists a homomorphism of ordered
groups $\varphi \colon H \to G$ such that $\varphi([e(n,k)]) = g(n,k)$
for all $n \ge 0$ and $0 \le k \le n$. 

Each function $g(n,k)$ is strictly positive, and hence non-zero, so it
follows from Proposition~\ref{prop:b1} that $\varphi$ is faithful. 
\end{proof} 

\noindent{\bf{Proof of Theorem~\ref{thm:GICAR}.}} By
Corollary~\ref{cor:initial} it suffices to find a unital embedding of
the GICAR-algebra $B$ into the AF-algebra $P = \bigotimes_{n=1}^\infty
M_2 \oplus M_3$. The ordered $K_0$-group of $P$ is (isomorphic to) a
dense subgroup $G$ of $C(X,\R)$ which contains the constant function
$1$ (as shown immediately after  Corollary~\ref{cor:initial}). 
The existence of a unital embedding of
the GICAR-algebra into the AF-algebra $P$ now follows from
Proposition~\ref{prop:d}. \hfill $\square$   

\section{Properties of initial objects} \label{sec:4}

We shall show in this last section that initial objects in the sense of 
this paper, although abundant, form at the same time  a rather  special 
class of \Cs s.

An element $g$ in an abelian group $G$ will be said to be \emph{infinitely
  divisible} if the set of natural numbers $n$ for which the equation
$nh=g$ has a solution $h \in G$ is unbounded. 

\begin{proposition} \label{prop:char1} If $B$ is an initial object,
  then $K_0(B)^+$ contains no non-zero infinitely divisible elements.
\end{proposition}

\begin{proof}
There exists  a unital \Cs{} $A$ of real rank zero and with no non-zero
finite-dimensional representations, such that no non-zero element in $K_0(A)$ is
infinitely divisible, and such that any non-zero projection
has a non-zero class in $K_0 (A)$. 
(For example, any irrational rotation
\Cs{}.) If $B$ is an initial object, then $B$ embeds into $A$, and
by choice of $A$ the corresponding ordered 
group homomorphism $K_0(B) \to K_0(A)$ takes any
non-zero positive element of $K_0(B)$ into a non-zero positive element
of  $K_0(A)$. Since the  image of an infinitely divisible element is again
infinitely divisible, no non-zero element of  $K_0(B)^+$ can be  
infinitely divisible.
\end{proof}

\begin{lemma} \label{lm:4}
Let $(G,G^+,u)$ be an ordered abelian group with order unit. Let $p_1,
p_2, \dots, p_n$ be distinct primes and suppose that $f_1, \dots, f_n$
are states on $(G,G^+,u)$ such that $f_j(G) = \Z[1/p_j]$ for $j=1,
\dots, n$. Then $f_1, \dots, f_n$ are affinely   independent.
\end{lemma}

\begin{proof}
The assertion  is proved by induction on $n$. It suffices to show that for
every natural number $n$, for every set of distinct primes $p_1,
\dots, p_n,q$, and for every set of states $f_1, \dots, f_n,f$ on
$(G,G^+,u)$, with $f_j(G) = \Z[1/p_j]$ and $f(G) = \Z[1/q]$ and with
$f_1, \dots, f_n$   affinely  independent,  $f$ is not an
affine combination of $f_1, \dots, f_n$.

Suppose, to reach a contradiction, that $f= \alpha_1f_1 + \cdots +
\alpha_nf_n$, with $\alpha_1, \dots, \alpha_n$  real numbers with
sum 1. If $n=1$, then $f=f_1$, which clearly is impossible. Consider the 
case $n \ge 2$. Since $f_1, \dots, f_n$ are assumed to be affinely independent,
there are $g_1, \dots, g_{n-1} \in G$ such that the vectors
$$x_j = (f_j(g_1), f_j(g_2), \dots, f_j(g_{n-1})) \in \Q^{n-1}, \qquad
j = 1, 2, \dots, n,$$
are affinely independent. The coefficients $\alpha_j$ above therefore
constitute the unique solution to the equations 
\begin{eqnarray*}
f_1(g_j)\alpha_1 + f_2(g_j)\alpha_2 + \cdots + f_n(g_j)\alpha_n &=& f(g_j), 
\qquad j=1,2, \dots, n-1,\\
\qquad \alpha_1 + \cdots + \alpha_n & = &1.
\end{eqnarray*}
As these $n$ equations in the $n$ unknowns $\alpha_1, \alpha_2, \dots,
\alpha_n$ are linearly independent, and all the
coefficients are rational, also $\alpha_1, \alpha_2, \dots, \alpha_n$
must be rational.

Denote by $\cQ'(q)$ the ring of all rational numbers with 
denominator (in reduced form) 
not divisible by $q$. Observe that $f_j(g) \in \cQ'(q)$ for all 
$j=1, \dots, n$ and for all $g \in G$. There is a natural number $k$ such that
$q^k\alpha_j \in \cQ'(q)$ for all $j=1, \dots, n$. Then
$$q^kf(g) = q^k\alpha_1f_1(g) + \cdots + q^k\alpha_nf_n(g) \in 
\cQ'(q),$$
for all $g \in G$. But this is impossible as, by hypothesis, 
$f(g) = 1/q^{k+1}$
for some $g \in G$.
\end{proof}

\begin{proposition} \label{prop:2}
Let $B$ be an initial object (in the sense of Definition~\ref{def:initial}),
and suppose that no quotient of $B$ has a 
minimal non-zero projection. Then the trace simplex $T(B)$ of $B$ is 
necessarily infinite-dimensional.
\end{proposition}

\noindent It follows in particular that any simple unital \Cs{} of real
rank zero, other than $\C$, which is an initial object has 
infinite-dimensional trace simplex. (Note for this that no matrix algebra
$M_n$ with $n \ge 2$ is an initial object.)

\begin{proof}
Any initial object embeds by definition into a large class of \Cs s that
includes exact \Cs s (such as for example any UHF-algebra),
and is therefore itself exact, being a sub-\Cs{} of an exact \Cs{} 
(see \cite[Proposition~7.1]{Kir:UHF}).
It follows (from \cite{BlaRor:extending} and \cite{Haa:quasi}, or from
\cite{HaaTho:traces}) that the canonical affine map from the trace
simplex $T(B)$ to the state space of 
$(K_0(B), K_0(B)^+,[1])$ is surjective.  
It is therefore sufficient to show that the latter space is 
infinite-dimensional. 
For each prime $p$ there is a unital embedding of $B$ into the UHF-algebra 
of type $p^\infty$, and hence a homomorphism of ordered groups
$f_p \colon K_0(B) \to
\Z[1/p]$ with $f_p([1]) = 1$.
Let us show that the homomorphisms $f_p$, when considered as states
(i.e., homomorphisms of ordered groups with order
unit from $(K_0 (A), [1])$ to $(\R, 1)$), are
affinely independent.

For each prime number $p$, the image of $f_p$ is a subgroup of
$\Z[1/p]$ which contains $1$, but  the only such subgroups are
$\Z[1/p]$ itself and the  subgroups  $p^{-k}\Z$ for some $k \ge 0$.
The latter cannot be the image of $f_p$ because the image of $B$ in our UHF-algebra, 
being isomorphic to a quotient of $B$, is assumed to have no minimal non-zero 
projections. (Indeed, if $\{p_n\}$ is a strictly decreasing sequence 
of projections in the sub-algebra of the UHF-algebra, and if $\tau$ is
the tracial state on the UHF-algebra, then  
$\{\tau(p_{n}-p_{n+1})\}$ is a sequence of strictly positive real numbers
which converges to $0$.)

Hence $f_p(K_0(B)) = \Z[1/p]$ for each prime $p$. It now follows from 
Lemma~\ref{lm:4} that the states $\{f_p : p\  {\text {prime}}\}$ 
are affinely independent. This shows that the state space of $(K_0(B),
K_0(B)^+,[1])$ is infinite-dimensional, as desired.
\end{proof} 

\noindent We end our paper by raising the following question:

\begin{problem} \label{pr}
Characterise initial objects (in the sense of Definition~\ref{def:initial})
among (simple) unital AF-algebras.
\end{problem}

\noindent We could of course extend the problem above to include all (simple)
real rank zero \Cs s, but we expect a nice(r) answer when we restrict
our attention to AF-algebras.
Propositions~\ref{prop:char1} and \ref{prop:2} give necessary, but not
sufficient, conditions for being an initial object. (A simple AF-algebra 
that satisfies the conditions of Propositions~\ref{prop:char1} and 
\ref{prop:2} can contain a unital simple sub-AF-algebra that does not 
satisfy the condition in Proposition~\ref{prop:2}, and hence is not an
initial object.)  

{\small{

\providecommand{\bysame}{\leavevmode\hbox to3em{\hrulefill}\thinspace}
\providecommand{\MR}{\relax\ifhmode\unskip\space\fi MR }
\providecommand{\MRhref}[2]{%
  \href{http://www.ams.org/mathscinet-getitem?mr=#1}{#2}
}
\providecommand{\href}[2]{#2}

}
\vspace{1cm}
Department of Mathematics\par
University of Toronto\par
Toronto, Canada~ M5S 3G3\par
e-mail address: elliott@math.toronto.edu

\vspace{.8cm} 
Department of Mathematics\par
University of Southern Denmark\par
5230 Odense M, Denmark\par
e-mail address: mikael@imada.sdu.dk
\end{document}